\documentclass{article}

\usepackage{amssymb}
\usepackage[cp1251]{inputenc}
\usepackage[english]{babel}
\usepackage{amsmath}
\usepackage{color}

\def\mapr#1{\smash{\mathop{\buildrel{#1}\over\longrightarrow}}}

\def\mapd#1{\Big\downarrow\rlap{$\vcenter{\hbox{$#1$}}$}}
\def\mapu#1{\Big\uparrow\rlap{$\vcenter{\hbox{$#1$}}$}}

\newtheorem{theorem}{Theorem}
\newtheorem{lemma}{Lemma}
\newtheorem{definition}{Definition}

\newtheorem{cor}{Corollary}
\def\qed{\hfill\vrule width2mm height2mm depth2mm}
\def\proof{{\bf Proof. }}
\def\f#1{{$#1$}}
\def\ff#1{{$$#1$$}\noindent}
\def\feq#1#2{\begin{equation}\label{#1}
    #2\end{equation}\noindent}
\def\({\left(}
\def\){\right)}

\def\C{{\bf C}}
\def\U{{\bf U}}

\def\id{{\hbox{\bf Id}}}
\def\1{{\bf 1}}

\title{Description of  the vector $G$-bundles over
$G$-spaces with quasi-free proper action of discrete group $G$}

\author{Mishchenko , A.S.\thanks{Partly supported by the
grant of RFFI No.08-01-00034-a,
NSh-1562.2008.1, Program 2.1.1/5031}\and Morales Mel\'endez, Quitzeh}
\begin{document}
\maketitle

\begin{abstract}
We give a description of the vector $G$-bundles over
$G$-spaces with qusi-free proper action of discrete group $G$ in terms
of the classifying space.
\end{abstract}

\section{The setting of the problem}

This problem naturally arises from the Conner-Floyd's description (\cite{Conner})
of the bordisms with the action of a group
$G$ using the so-called fix-point construction.
This construction reduces the problem of describing the bordisms to two simpler
problems: a) description of the fixed-point set (or, more generally, the stationary point set),
which happens to be a submanifold attached with the structure of its normal bundle
and the action of the same group $G$, however, this action could have stationary points
of lower rank; b) description of the bordisms of lower rank with an action of the group $G$.
We assume that the group  $G$ is discrete.

Lets  $\xi$ be an $G$-equivariant vector bundle with base $M$.
\feq{1}{
\begin{array}{c}
          \xi \\
          \mapd{} \\
          M \\
        \end{array}
}Lets $H<G$ be a normal finite subgroup. Assume that the action of the group
$G$  over the base $M$ reduces to the factor group $G_{0}=G/H$:
\feq{2}{
\begin{array}{ccc}
  G\times M &\mapr{} &M \\
  \mapd{} &&\parallel\\
  G_{0}\times M&\mapr{}& M\\
\end{array}}suppose, additionally,
that the action $G_{0}\times M\mapr{} M$ is free and there is no more fixed points
of the action of the group $H$ in the total space of the bundle $\xi$.

So, we have the following commutative diagram

\feq{3}{
\begin{array}{ccc}
  G\times\xi &\mapr{}&\xi\\
  \mapd{} &&\mapd{}\\
  G_{0}\times M&\mapr{}&M \\
\end{array}
}

\begin{definition}
 As in \cite[p. 210]{Levine}, we shall say that the described
 action of the group $G$ is \textit{quasi-free}
over the base with normal \textit{stationary} subgroup $H$.
\end{definition}

Reducing the action to the subgroup $H$, we obtain the simpler diagram:
\feq{4}{
\begin{array}{ccc}
  H\times\xi &\mapr{}&\xi\\
  \mapd{} &&\mapd{}\\
   M&=&M \\
\end{array}
}

Following \cite{Atiyah}, let $\rho_{k}: H\mapr{}\U(V_{k})$ be the series of all the irreducible
(unitary) representation of the finite group $H$. Then the $H$-bundle $\xi$
can be presented as the finite direct sum:
\feq{5}{
\xi\approx \bigoplus_{k}\(\xi_{k}\bigotimes V_{k}\),
}
where the action of the group  $H$ over the bundles $\xi_{k}$ is trivial,  $V_{k}$
denotes the trivial bundle with fiber $V_{k}$ and with  fiberwise
action of the group $H$, defined using the linear representation $\rho_{k}$.

\begin{lemma} The group $G$ acts on every term of the sum (\ref{5})
separately.
\end{lemma}

\proof Consider now the action of the group $G$ over the total space
of the bundle $\xi$. Fix a point $x\in M$. The action of the element
$g\in G$ is fiberwise, and maps the fiber $\xi_{x}$ to the fiber
$\xi_{gx}$: \ff{ \Phi(x,g):\xi_{x}\mapr{}\xi_{gx}. }

 Also, for a par of elements  $g_{1}, g_{2}\in G$ we have:

\feq{5.1}{
\Phi(x, g_{1}g_{2})=\Phi\(g_{2}x, g_{1}\)\circ\Phi\(x, g_{2}\),
}

\ff{
\begin{array}{ccccc}
 \Phi(x, g_{1}g_{2}):     \xi_{x}& \mapr{\Phi\(x, g_{2}\)}&\xi_{g_{2}x}&
      \mapr{\Phi\(g_{2}x, g_{1}\)} & \xi_{g_{1}g_{2}x}\\
\end{array}
}

In particular, if $g_{2}=h\in H<G$, then $g_{2}x=hx=x$. So,
\ff{
\begin{array}{ccccc}
 \Phi(x, gh):     \xi_{x}& \mapr{\Phi\(x, h\)}&\xi_{x}&
      \mapr{\Phi\(x, g\)} & \xi_{gx}\\
\end{array}
} Analogously, if $g_{1}=h\in H<G$, then $g_{1}gx=hgx=gx$. So \ff{
\begin{array}{ccccc}
 \Phi(x, hg):     \xi_{x}&
      \mapr{\Phi\(x, g\)} & \xi_{gx}& \mapr{\Phi\(gx, h\)}\xi_{gx}\\
\end{array}
}
According to \cite{Atiyah} the operator $\Phi\(x, h\)$ does not depends on the point $x\in M$,

\ff{
\Phi(x, h)=\Psi(h):\bigoplus_{k}\(\xi_{k,x}\bigotimes V_{k}\)\mapr{}
\bigoplus_{k}\(\xi_{k,x}\bigotimes V_{k}\),
}here, since the action of the group  $H$ is given over every space
$V_{k}$ using pairwise different irreducible representations $\rho_{k}$,
we have
\ff{
\Psi(h)=\bigoplus_{k}\(\id\bigotimes \rho_{k}(h)\).
}

In this way, we obtain the following relation:
\feq{6}{
\Phi(x, gh)=\Phi(x, g)\circ\Psi(h)=
\Phi(x, ghg^{-1}g)=\Psi(ghg^{-1})\circ\Phi(x, g).
}

Lets write the operator  $\Phi(x, g)$ using matrices to decompose the space
 $\xi_{x}$ as the direct sum
\ff{
\xi_{x}=\bigoplus_{k}\(\xi_{k,x}\bigotimes V_{k}\):
}
\feq{6.1}{
\Phi(x, g)=\(
\begin{array}{cccc}
   \Phi(x, g)_{1,1}& \cdots & \Phi(x, g)_{k,1} & \cdots \\
  \vdots & \ddots & \vdots &  \\
  \Phi(x, g)_{1,k} & \cdots & \Phi(x, g)_{k,k} & \cdots\\
  \vdots &  & \vdots & \ddots \\
\end{array}
\)
}

If $k\neq l$ then $\Phi(x, g)_{k,l}=0$, i.e. the matrix $\Phi(x, g)$ its diagonal,
\ff{
\Phi(x, g)=\bigoplus_{k}\Phi(x, g)_{k,k}:
\bigoplus_{k}\(\xi_{k,x}\bigotimes V_{k}\)\mapr{}
\bigoplus_{k}\(\xi_{k,gx}\bigotimes V_{k}\),
}
\ff{
\Phi(x, g)_{k,k}:
\(\xi_{k,x}\bigotimes V_{k}\)\mapr{}
\(\xi_{k,gx}\bigotimes V_{k}\),
}
as it was required to prove. \qed

\section{Description of the particular case $\xi=\xi_{0}\bigotimes V$}
Here we will consider the particular case of a $G$-vector bundle $\xi=\xi_{0}\otimes V$ with base $M$.
\ff{
\begin{array}{c}
          \xi \\
          \mapd{} \\
          M \\
        \end{array}
}
where  the action of the group $G$ is quasi-free over the base with finite normal stationary subgroup  $H<G$.

We will assume that the group $H$ acts trivially over the bundle $\xi_0$.
By $V$ we denote the trivial bundle with fiber $V$ and with fiberwise action
of the group $H$ given by an irreducible linear representation $\rho$.

\begin{definition} A \textit{canonical model} for the fiber in a $G$-bundle
$\xi=\xi_0\bigotimes V$ with fiber $F\otimes V$ is the product
$G_0\times \(F\otimes V\)$ with an action of the group $G$ \ff{
\begin{array}{ccc}
  G\times \(G_0\times \(F\otimes V\)\) &\mapr{\phi}& G_0\times \(F\otimes V\)\\
  \mapd{} &&\mapd{}\\
  G\times G_{0}&\mapr{\mu}&G_{0} \\
\end{array}
}where
$\mu$ denotes the natural left action of $G$ on its quotient $G_0$, and
\ff{
\phi([g],g_{1}):[g]\times \(F\otimes V\) \to [g_1g]\times \(F\otimes V\)
}is
given by the formula
\feq{15.1a}{
\begin{array}{ll}
\phi([g],g_{1})
=
\id\otimes
\rho(u(g_{1}g)u^{-1}(g)).
    \end{array}
}where
\ff{
u:G\mapr{}H
}
is a homomorphism of right $H$-modules by multiplication,
i.e.
\ff{
u(gh)=u(g)h,\quad  u(1)=1, \quad g\in G, h\in H.
}
\end{definition}

\begin{lemma} The definition (\ref{15.1a}) of the action of  \f{G} is
well-defined.
\end{lemma}

\proof It is enough to prove that
that  a) the formula  (\ref{15.1a}) defines an action, i.e.
\ff{
\begin{array}{ll}
\phi([g],g_2g_{1})
=\phi([g_1g],g_2)\circ\phi([g],g_{1}),
    \end{array}
}and b) that the formula (\ref{15.1a}) does not depends on the chosen representative $gh\in [g]$:
\ff{\id\otimes \rho(u(g_{1}g)u^{-1}(g))=\id\otimes \rho(u(g_{1}gh)u^{-1}(gh))}for
 every $g\in G$ and $h\in H$.

In fact,
\ff{
\begin{array}{ll}
\phi([g],g_2g_{1})
=
\id\otimes \rho(u(g_2g_{1}g)u^{-1}(g))=\\\\
\id\otimes \rho(u(g_2g_{1}g)u(g_{1}g)u^{-1}(g_{1}g)u^{-1}(g))=\\
=\id\otimes \rho(u(g_2g_{1}g)u(g_{1}g))\circ\id\otimes \rho(u^{-1}(g_{1}g)u^{-1}(g))=\\\\
=\phi([g_1g],g_2)\circ\phi([g],g_{1}),
    \end{array}
}what proves a), and,  recalling the equation $u(gh)=u(g)h$ for every $g\in G$ and $h\in H$, it is clear that
\ff{u(g_{1}gh)u^{-1}(gh)=u(g_{1}g)hh^{-1}u^{-1}(g)=u(g_{1}g)u^{-1}(g),}
 which is a sufficient condition for b) to be true. \qed

\pagebreak
As it is well known, for the actions we are studying, we can always consider over the base $M$ an atlas of
equivariant charts $\{O_{\alpha}\}$,
\ff{
M=\bigcup_{\alpha}O_{\alpha},
}
\ff{
[g]O_{\alpha}=O_{\alpha}, \qquad\forall [g]\in G_{0}.
}
If the atlas is fine enough, then every chart can be presented as a
disjoint union of its subcharts:
\ff{
O_{\alpha}= \bigsqcup_{[g]\in G_{0}}[g]U_{\alpha}\approx U_{\alpha}\times G_{0},
}i.e. \f{[g]U_{\alpha}\cap [g']U_{\alpha}=\emptyset} if \f{ [g]\neq [g']}, and
when $\alpha\neq\beta$,  if
$U_{\alpha}\cap [g_{\alpha\beta}]U_{\beta}\neq \emptyset$,
then the element $g_{\alpha\beta}$ is the only one for which that intersection is
non-empty,
i.e. if $[g]\neq[g_{\alpha\beta}]$,
then $U_{\alpha}\cap[g]U_{\beta}=\emptyset$,
i.e.
\ff{
O_{\alpha}\cap O_{\beta} \approx \(U_\alpha\cap [g_{\alpha\beta}]U_\beta\)\times G_0,
}for every $\alpha, \beta$. We use these facts and notations to formulate the next theorem.

\begin{theorem}
The bundle  $\xi= \xi_{0}\bigotimes V$ is locally homeomorphic
to the cartesian product of some chart $U_\alpha$ by the canonical
model.  More precisely, for a fine enough atlas, there exist
$G$-equivariant trivializations
\feq{14.a}{\psi_\alpha:O_\alpha\times \(F\otimes
V\)\to\xi|_{O_\alpha}}where \ff{O_\alpha \times \(F\otimes
V\)\approx U_\alpha \times \(G_0\times \(F\otimes V\)\)}and the
diagram \feq{trivial}{
\begin{array}{ccc}\xi|_{O_\alpha}&\mapr{g}&\xi|_{O_\alpha}\\
                               \mapu{\psi_\alpha} &&\mapu{\psi_\alpha}\\
                               U_\alpha\times \(G_0\times \(F\otimes V\)\)
                               &\mapr{\id\times \phi(g)}&U_\alpha\times \(G_0\times \(F\otimes V\)\)\\
\end{array}
}is
commutative where $g\in G$, \f{\id:U_\alpha\to U_\alpha,} and $\phi(g)$ denotes
the canonical action.
\end{theorem}

\proof Using an atlas as in the remarks at the beginning of the
theorem, we shall construct the trivialization (\ref{14.a}) starting
from an arbitrary trivialization \ff{\psi_\alpha:U_\alpha\times
\(F\otimes V\)\to\xi|_{U_\alpha}}in such a way, that the diagram
\ff{
\begin{array}{ccc}\xi|_{U_\alpha}&\mapr{g}&\xi|_{[g]U_\alpha}\\
                               \mapu{\psi_\alpha} &&\mapu{\psi_\alpha}\\
                               U_\alpha\times \(F\otimes V\)
                               &\mapr{}&[g]U_\alpha\times  \(F\otimes V\)\\
\end{array}
}commutes for every $g\in [g]$, where the left and upper arrows
are given and we have to construct the down and right arrows.

From such a construction,
the equivariance will follow automatically  and the proof of the theorem reduces to show that the
constructed down arrow coincides with that on (\ref{trivial}).

Evidently, for a given trivialization \f{\psi_\alpha:U_\alpha\times \(F\otimes V\)\to\xi|_{U_\alpha}},
there are several ways to define a trivialization
\f{\psi_\alpha:[g]U_\alpha\times \(F\otimes V\)\to\xi|_{[g]U_\alpha}}, since there are several
elements $g\in G$ sending $\xi|_{U_\alpha}$ to $\xi|_{[g]U_\alpha}$.

Thus, consider a set-theoretic cross-section
\ff{
p':G_{0}\mapr{}G,
}
to the projection $p$ in the exact sequence of groups
\ff{
\1\mapr{}H\mapr{}G\mapr{p}G_{0},
}
\ff{
p\circ p' = \id : G_{0}\mapr{p'}G\mapr{p}G_{0}.
}
Put
\ff{
g'=p'\circ p: G\mapr{}G.
}
Without loss of generality, we can take  $g'(1)=1$.

In this case
\ff{
g'(g)=gu^{-1}(g),
}where
\ff{
u:G\mapr{}H
}
is a homomorphism of right $H$-modules by multiplication,
i.e.
\ff{
u(gh)=u(g)h, \quad g\in G, h\in H.
}
In particular, this means that
\ff{
g'(gh)=g'(g), \quad h\in H.
}

Lets
\ff{
\tilde\psi_{\alpha}:U_{\alpha}\times F
\mapr{} \xi_0|_{U_{\alpha}}
}be
some trivialization. We define the trivialization $\psi_{\alpha}$ in (\ref{14.a})
by the rule: if \f{[g]x_{\alpha}\in [g]U_{\alpha}}, i.e.
\f{x_{\alpha}\in U_{\alpha}},
then, the map

\ff{
\psi_{\alpha}([g]x_{\alpha}):
[g]x_{\alpha}\times \(F\otimes V\)
\mapr{} \xi_{[g]x_{\alpha}}\otimes V
}is given by the formula
\feq{15a}{\begin{array}{ll}
\psi_{\alpha}([g]x_{\alpha})&=
\Phi(x_{\alpha}, g'(g))\circ \(\tilde\psi_{\alpha}(x_{\alpha})\otimes \id\)=\\
&=\Phi(x_{\alpha}, gu^{-1}(g))\circ
\(\tilde\psi_{\alpha}(x_{\alpha})\otimes \id\).\\
\end{array}
}where, from the first equality, it is clear that the definition does not depend
on the representative $g\in [g]$.

In particular, for $[g]=1$, we recover the initial trivialization
\ff{
\psi_{\alpha}(x_{\alpha})=\tilde\psi_{\alpha}(x_{\alpha})\otimes \id
}
since \f{\Phi(x, g'(1))=\Phi(x, 1)= 1}.

Using this trivialization the action of the group $G$ can be carried to
the cartesian product \f{O_{\alpha}\times \(F\otimes V\)}:
\ff{
\Phi_{\alpha}(g):O_{\alpha}\times \(F\otimes V\)
\mapr{}O_{\alpha}\times \(F\otimes V\).
}
Lets \f{x_{\alpha}\in U_{\alpha}}, \f{g\in G}, then
\ff{
\Phi_{\alpha}([g]x_{\alpha},g_{1}):[g]x_{\alpha}
\times\(F\otimes V\)\mapr{}[g_{1}g]x_{\alpha}\times\(F\otimes V\)
}
is given by the formula
\ff{
\Phi_{\alpha}([g]x_{\alpha},g_{1})=
\(\psi_{\alpha}([g_{1}g]x_{\alpha})\)^{-1}
\Phi([g]x_{\alpha},g_{1})\psi_{\alpha}([g]x_{\alpha}).
}

Applying  (\ref{15a}), we obtain
\ff{\begin{array}{ll}
      \Phi_{\alpha}([g]x_{\alpha},g_{1})= &
      \(\Phi(x_{\alpha}, g_{1}gu^{-1}(g_{1}g))\circ
      \(\tilde\psi_{\alpha}(x_{\alpha})\otimes \id\)\)^{-1}
\circ\\&
\circ\Phi([g]x_{\alpha},g_{1})\circ
\Phi(x_{\alpha}, gu^{-1}(g))\circ \(\tilde\psi_{\alpha}(x_{\alpha})\otimes \id\) =\\ \\
&=\(\tilde\psi_{\alpha}(x_{\alpha})\otimes \id\)^{-1}\circ \\
&\circ
\Phi(x_{\alpha}, g_{1}gu^{-1}(g_{1}g))^{-1}\circ
\Phi([g]x_{\alpha},g_{1})\circ\Phi(x_{\alpha}, gu^{-1}(g))
\circ\\
&\circ\(\tilde\psi_{\alpha}(x_{\alpha})\otimes \id\)=\\\\
&=\(\tilde\psi_{\alpha}(x_{\alpha})\otimes \id\)^{-1}\circ \\
&\circ
\Phi(x_{\alpha}, u^{-1}(g_{1}g))^{-1}\circ
\Phi(x_{\alpha}, g_{1}g)^{-1}\circ\Phi([g]x_{\alpha},g_{1})\circ\\
&\circ
\Phi(x_{\alpha}, g)\circ\Phi(x_{\alpha}, u^{-1}(g))
\circ\\
&\circ\(\tilde\psi_{\alpha}(x_{\alpha})\otimes \id\)=\\\\
&=\(\tilde\psi_{\alpha}(x_{\alpha})\otimes \id\)^{-1}\circ \\
&\circ
\Phi(x_{\alpha}, u^{-1}(g_{1}g))^{-1}\circ
\Phi(x_{\alpha}, u^{-1}(g))
\circ\\
&\circ\(\tilde\psi_{\alpha}(x_{\alpha})\otimes \id\);\\\\
\end{array}
}

\ff{\begin{array}{ll}
\Phi_{\alpha}([g]x_{\alpha},g_{1})
=&\(\tilde\psi_{\alpha}(x_{\alpha})\otimes \id\)^{-1}\circ \\
&\circ
\(\id\otimes\rho(u(g_{1}g))\)\circ
\(\id\otimes\rho(u^{-1}(g))\)
\circ\phantom{aaaaaaaaaaaaaaaaa}\\
&\circ\(\tilde\psi_{\alpha}(x_{\alpha})\otimes \id\)=\\\\
&=\(\tilde\psi_{\alpha}(x_{\alpha})\otimes \id\)^{-1}\circ \\
&\circ
\(\id\otimes\(
\rho(u(g_{1}g)u^{-1}(g))
\)\)
\circ\\
&\circ\(\tilde\psi_{\alpha}(x_{\alpha})\otimes \id\)=\\\\
&=
\id\otimes
\rho(u(g_{1}g)u^{-1}(g))
.
    \end{array}
}

The operator
\ff{
\begin{array}{ll}
\Phi_{\alpha}([g]x_{\alpha},g_{1})
=
\id\otimes
\rho(u(g_{1}g)u^{-1}(g))=\phi(g_{1}, [g]).
    \end{array}
}
does not depend on the point $x_{\alpha}\in U_{\alpha}$. So, the theorem is proved.\qed

By \f{\mathrm{Aut}_G\(G_0\times \(F\otimes V\)\)} we denote the group of equivariant automorphisms
of the space $G_0\times \(F\otimes V\)$ as a vector $G$-bundle with base $G_0$, fiber
$F\otimes V $ and canonical action of the group $G$.

\begin{cor} The transition functions on the intersection \ff{ O_\alpha \cap
O_\beta\approx \(U_\alpha\cap [g_{\alpha\beta}]U_\beta\)\times G_0,}
i.e. the homomorphisms \f{\Psi_{\alpha\beta}} on the diagram
\feq{15.2a}{
\begin{array}{ccc}
  \(U_\alpha\cap  [g_{\alpha\beta}]U_\beta\)\times \(G_0\times \(F\otimes V\)\)
  &\mapr{\Psi_{\alpha\beta}}& \(U_\alpha\cap  [g_{\alpha\beta}]U_\beta\)\times \(G_0\times \(F\otimes V\)\)\\
  \mapd{} &&\mapd{}\\
   \(U_\alpha\cap  [g_{\alpha\beta}]U_\beta\)\times G_0
  &\mapr{\id}& \(U_\alpha\cap [g_{\alpha\beta}]U_\beta\)\times G_0\\
\end{array}
}are equivariant with respect to the canonical action of the group
$G$ over the product of the base by the canonical model, i.e.
\ff{\Psi_{\alpha\beta}(x)\circ
\phi(g_1,[g])=\phi(g_1,[g])\circ\Psi_{\alpha\beta}(x)}for every
$x\in U_\alpha\cap  [g_{\alpha\beta}]U_\beta,\;g_1\in G,\; [g]\in
G_0$, In other words, \ff{\Psi_{\alpha\beta}(x)\in
\mathrm{Aut}_G\(G_0\times \(F\otimes V\)\).}
\end{cor}

Now we give a more accurate description of the group
\f{\mathrm{Aut}_G\(G_0\times \(F\otimes V\)\)}.
By definition, an element of the group \f{\mathrm{Aut}_G\(G_0\times \(F\otimes V\)\)}
is an equivariant mapping $\mathbf{A}^{a}$, such that the pair $(\mathbf{A}^{a},a)$
defines a commutative diagram
\ff{
\begin{array}{ccc}
   \(G_0\times \(F\otimes V\)\) &\mapr{\mathbf{A}^{a}}& G_0\times \(F\otimes V\)\\
  \mapd{} &&\mapd{}\\
   G_{0}&\mapr{a}&G_{0}, \\
\end{array}
}
which commutes with the canonical action, i.e.
the map \f{a\in \mathrm{Aut}_G(G_0)} satisfies the condition
\ff{a\in \mathrm{Aut}_G(G_0)\approx G_0,\quad a[g]=[ga],\;[g]\in G_0,}
and the mapping   $\mathbf{A}^{a}=(A^{a}[g])_{[g]\in G_0}$,

\ff{
A^{a}[g]: [g]\times(F\otimes V)\to [ga]\times (F\otimes V)
}
satisfies a commutation condition with respect to the action of the group \f{G}:
\ff{
\begin{array}{ccc}
  [g]\times(F\otimes V) & \mapr{A^{a}[g]} & [ga]\times(F\otimes V) \\
  \mapd{\phi(g_{1},[g])} &  & \mapd{\phi(g_{1},[ga])} \\
  {[g_{1}g]\times(F\otimes V)} & \mapr{A^{a}[g_{1}g]} & [g_{1}ga]\times(F\otimes V)
\end{array} \quad ,
}

\feq{16.e}{\phi(g_1,[ga])\circ A^{a}[g]=A^{a}[g_{1}g]\circ\phi(g_1,[g])}
 i.e.
 \feq{16.b}{
 (\id\otimes \rho(u(g_{1}ga)u^{-1}(ga)))A^{a}[g]=
 A^{a}[g_1g](\id\otimes\rho(u(g_{1}g)u^{-1}(g)))
 }
where  $ [g]\in G_0, \quad g_1\in G$.

\begin{lemma} One has an exact sequence of groups
\feq{suc}{\1\to GL(F)
\mapr{}\mathrm{Aut}_G\(G_0\times \(F\otimes V\)\)
\mapr{} G_0\to \1.}
\end{lemma}

\proof To define a projection
\ff{pr:\mathrm{Aut}_G\(G_0\times \(F\otimes V\)\) \mapr{} G_0}we
send the fiberwise map
\ff{\mathbf{A}^{a}:G_0\times \(F\otimes V\) \mapr{} G_0\times \(F\otimes V\)}to
its restriction over the base
\f{a:G_0\to G_0}, i.e. $a\in \mathrm{Aut}_G(G_0)\approx G_0$. So, this is a well-defined homomorphism.

We need to show that $pr$ is an epimorphism and that its kernel is isomorphic to $GL(F)$.
Lets calculate the kernel.

For $[a]=[1]$ we have
 \feq{17.a}{
 (\id\otimes \rho(u(g_{1}g)u^{-1}(g)))A^{1}[g]=
 A^{1}[g_1g](\id\otimes\rho(u(g_{1}g)u^{-1}(g)))
 }

In the case  $g_{1}=h\in H$, we obtain

 \ff{
 (\id\otimes \rho(u(hg)u^{-1}(g)))A^{1}[g]=
 A^{1}[g](\id\otimes\rho(u(hg)u^{-1}(g)))
 }
Since the representation \f{\rho} is irreducible,
by Schur's lemma, we have
\ff{
A^{1}[g] = B^{1}[g]\otimes \id.}

On the other side, assuming in (\ref{17.a}) that \f{g=1}, we have

 \ff{
 (\id\otimes \rho(u(g)))A^{1}[1]=
 A^{1}[g](\id\otimes\rho(u(g))),
 }
 i.e.

  \ff{
 (\id\otimes \rho(u(g)))(B^{1}[1]\otimes \id)=
 (B^{1}[g]\otimes \id)(\id\otimes\rho(u(g))),
 }
 or
   \ff{
 (B^{1}[g]\otimes \id)=
 (B^{1}[1]\otimes \id).
 }

 So, the kernel \f{\ker pr} is isomorphic to the group \f{GL(F)}.

 In the generic case, i.e. \f{[a]\neq 1},
 we can compute the operator \f{A^{a}[g]} in terms of its value at the identity \f{A^{a}[1]}
 from the formula (\ref{16.b}):
assuming  $g=1$, we obtain (changing $g_{1}$ by $g$):

  \feq{16.c}{
 (\id\otimes \rho(u(ga)u^{-1}(a)))A^{a}[1]=
 A^{a}[g](\id\otimes\rho(u(g))),
 }

 i.e.

  \feq{16.d}{
A^{a}[g]=
(\id\otimes \rho(u(ga)u^{-1}(a)))A^{a}[1](\id\otimes\rho(u^{-1}(g))),
 }

 Therefore, the operator is completely defined by its value
\ff{
A^{a}[1]: [1]\times(F\otimes V)\to [a]\times (F\otimes V)
}
at the identity $g=1$.

Now we describe the operator $A^{a}[1]$ in terms of the representation
$\rho$ and its properties.

We have a commutation rule with respect to the action of the subgroup \f{H}:
\ff{
\begin{array}{ccc}
  [1]\times(F\otimes V) & \mapr{A^{a}[1]} & [a]\times(F\otimes V) \\
  \mapd{\phi(h,[1])} &  & \mapd{\phi(h,[a])} \\
  {[1]\times(F\otimes V)} & \mapr{A^{a}[1]} & [a]\times(F\otimes V)
\end{array} \quad ,
}
Equivalently

\ff{
A^{a}[1]\circ\phi(h,[1])=\phi(h,[a])\circ A^{a}[1],
}

i.e.

\ff{
A^{a}[1]\circ(\id\otimes \rho(h))=(\id\otimes\rho(g'^{-1}(a)hg'(a)))\circ A^{a}[1],
}

i.e.

\ff{
A^{a}[1]\circ(\id\otimes \rho(h))=(\id\otimes\rho_{g'(a)}(h))\circ A^{a}[1].
}

The last equation means that the operator should $A^{a}[1]$ permute these
representations, or equivalently, such an operator exists only when the
representations $\rho$ and $\rho_{g'(a)}$ are equivalent. Recalling
the commutation rule (\ref{6}), we see that this is the case we are been
considering.

Thus, if the representations  $\rho$ and $\rho_{g}$ are equivalent,
 we have an (inverse) splitting operator \f{C(g)},
satisfying the equation
\feq{29-a}{
\rho_{g}(h)=\rho\(g^{-1}hg\)=C(g)\rho(h)C^{-1}(g).
}for every $g\in G$.
The operator  \f{C(g)} is defined up to multiplication by a scalar operator
\f{\mu_{g}\in\SS^{1}\subset\C^{1}}.

So

\ff{
A^{a}[1]\circ(\id\otimes \rho(h))=(\id\otimes C(g'(a))\circ\rho(h)\circ C^{-1}(g'(a)))\circ A^{a}[1],
}
or

\ff{
(\id\otimes C^{-1}(g'(a)))\circ A^{a}[1]\circ(\id\otimes \rho(h))=(\id\otimes\rho(h))\circ (\id\otimes C^{-1}(g'(a)))\circ A^{a}[1],
}

Then, by the Schur's lemma,

\ff{
(\id\otimes C^{-1}(g'(a)))\circ A^{a}[1]=B^{a}[1]\otimes\id,
}
i.e.
\ff{
A^{a}[1]=B^{a}[1]\otimes C(g'(a)),
}

Using the formula (\ref{16.d}), we obtain
  \ff{
A^{a}[g]=
(\id\otimes \rho(u(ga)u^{-1}(a)))(B^{a}[1]\otimes C(g'(a)))(\id\otimes\rho(u^{-1}(g))),
 }

i.e.

  \feq{16.f}{
A^{a}[g]=
B^{a}[1]\otimes (\rho(u(ga)u^{-1}(a))\circ C(g'(a))\circ\rho(u^{-1}(g))).
 }

This means, that by defining the matrix $B^{a}[1]$, it is possible to obtain all the operators
$A^{a}[g]$ satisfying the equation (\ref{16.d}).

It remains to verify the commutation rule (\ref{16.b}), i.e. in the formula

 \ff{
 (\id\otimes \rho(u(g_{1}ga)u^{-1}(ga)))A^{a}[g]=
 A^{a}[g_1g](\id\otimes\rho(u(g_{1}g)u^{-1}(g)))
 }
we substitute the expression (\ref{16.f}):
\ff{
\begin{array}{c}
  (\id\otimes \rho(u(g_{1}ga)u^{-1}(ga)))\circ(
  B^{a}[1]\otimes (\rho(u(ga)u^{-1}(a))\circ C(g'(a))\circ\rho(u^{-1}(g)))) =\\\\
  =(B^{a}[1]\otimes (\rho(u(g_{1}ga)u^{-1}(a))\circ C(g'(a))\circ\rho(u^{-1}(g_{1}g))))
  \circ(\id\otimes\rho(u(g_{1}g)u^{-1}(g)))
\end{array}
}

that is
\ff{
\begin{array}{c}
  B^{a}[1]\otimes \rho(u(g_{1}ga)u^{-1}(ga)))\circ
  (\rho(u(ga)u^{-1}(a))\circ C(g'(a))\circ\rho(u^{-1}(g)))) =\\\\
  =B^{a}[1]\otimes (\rho(u(g_{1}ga)u^{-1}(a))\circ C(g'(a))\circ\rho(u^{-1}(g_{1}g))))
  \circ(\rho(u(g_{1}g)u^{-1}(g)))
\end{array}
}

Note that this identity does not depend on the particular matrix $B^{a}[1]$,
thus, this means that we only need to verify the identity for arbitrary $a,g$ and $g_{1}$:

\ff{
\begin{array}{c}
  \rho(u(g_{1}ga)u^{-1}(ga)))\circ
  (\rho(u(ga)u^{-1}(a))\circ C(g'(a))\circ\rho(u^{-1}(g)))) =\\\\
  =(\rho(u(g_{1}ga)u^{-1}(a))\circ C(g'(a))\circ\rho(u^{-1}(g_{1}g))))
  \circ(\rho(u(g_{1}g)u^{-1}(g))),
\end{array}
}

which is obvious, after the natural simplifications

\ff{
\begin{array}{c}
  \rho(u(g_{1}ga)u^{-1}(a))\circ C(g'(a))\circ\rho(u^{-1}(g)))) =\\\\
  =(\rho(u(g_{1}ga)u^{-1}(a))\circ C(g'(a))\circ\rho(u^{-1}(g))),
\end{array}
}

So, it follows, that for every element $[a]\in G_{0}$ there exist an element
$(A^{a}, a)\in \mathrm{Aut}_G\(G_0\times \(F\otimes V\)\)$.
This means that the homomorphism
\ff{\mathrm{Aut}_G\(G_0\times \(F\otimes V\)\)
\mapr{pr} G_0}
is in fact an epimorphism, and the lemma is proved. \qed

 It is clear that there is an  equivalence between $G$-vector bundles with fiber
 $G_0\times \(F\otimes V\) $ over a (compact)
 base $X$, where $G$ acts trivially over the base and canonically over the fiber,
 and homotopy classes of mappings from $X$ to the space
 $B\mathrm{Aut}_G\(G_0\times \(F\otimes V\)\)$.

Lets denote by $\mathrm{Vect}_G(M,\rho)$ the category of $G$-equivariant vector
bundles $\xi=\xi_{0}\otimes V$ with base $M$, where  the action of the group $G$
is quasi-free over the base with finite normal stationary subgroup  $H<G$,
 the group $H$ acts trivially over the bundle $\xi_0$ and $V$ denotes
 the trivial bundle with fiber $V$ and with fiberwise action
of the group $H$ given by an irreducible linear representation $\rho$.
Here we need to require for the representations
\f{\rho_g(h)=\rho(g^{-1}hg)} to be equivalent for every $g\in G$, in the other case,
in view of the commutation rule, this category may be void.

This is a category because, in fact,
we are just taking vector bundles over the space
$M$, then applying tensor product by the fixed bundle $V$ and defining some action
of the group $G$ over the resulting spaces. The inclusion
$GL(F)\hookrightarrow \mathrm{Aut}_G\(G_0\times \(F\otimes V\)\)$
from lemma 2 ensures that the identities are included.

Denote by $\mathrm{Bundle}(X,L)$ the category of principal $L$-bundles over the
base $X$.

\begin{theorem} There is a monomorphism
\feq{categ}{\mathrm{Vect}_G(M,\rho)\longrightarrow
\mathrm{Bundle}(M/G_0,\mathrm{Aut}_G\(G_0\times \(F\otimes V\)\)). }
\end{theorem}

\proof By corollary 3, every element $\xi\in
\mathrm{Vect}_G(M,\rho)$ is defined by transition functions \ff{
\Psi_{\alpha\beta}:\; \(U_\alpha\cap [g_{\alpha\beta}]U_\beta\)\to
\mathrm{Aut}_G\(G_0\times \(F\otimes V\)\)}where by construction,
when $[g]\neq[g_{\alpha\beta}]$, we have
$U_{\alpha}\cap[g]U_{\beta}=\emptyset$ and if $[g]\neq 1$, then
$U_{\alpha}\cap[g]U_{\alpha}=\emptyset$ and
 $U_{\beta}\cap[g]U_{\beta}=\emptyset$. This means that
 the sets $U_{\alpha}$ and $U_{\beta}$  project homeomorphically to open sets under
 the natural projection $M\to M/G_0$. So, these transition
 functions are well-defined over an atlas of the quotient space $M/G_0$ and
 they form a $G$-bundle with fiber $G_0\times \(F\otimes V\)$
 over  this quotient space.

 By the same arguments, it is obvious that every $G$-equivariant map
 \feq{tr1}{h_\alpha:O_\alpha\times \(F\otimes V\)\to O_\alpha\times \(F\otimes V\)}can
be interpreted as a map
\feq{tr2}{h_\alpha:U_\alpha\times\(G_0\times \(F\otimes V\)\)\to U_\alpha\times\(G_0\times \(F\otimes V\)\)}by
means of the homeomorphism \f{O_\alpha\approx U_\alpha\times G_0}, where the set $U_\alpha$
can be thought as an open set of the space $M/G_0$. Equivalently,
\feq{tr3}{h_\alpha:U_\alpha \to \mathrm{Aut}_G\(G_0\times \(F\otimes V\)\)}where
$U_\alpha$ is homeomorphic to an open set of the space $M/G_0$.
Therefore, the map
(\ref{categ}) is well defined.

Conversely, if we start from mappings of the form
(\ref{tr3}) where the sets $U_\alpha$ are open in $M/G_0$, by refining the
atlas, if it is necessary, we can always think that the inverse image of the open sets $U_\alpha$
under the quotient map $M\to M/G_0$ are homeomorphic to the product $U_\alpha\times G_0$
and then obtain mappings of the form (\ref{tr1}). Therefore, the map
(\ref{categ}) is a monomorphism. \qed

Of course, the map (\ref{categ}) its not in general an epimorphism, since, when we define the
category $\mathrm{Vect}_G(M,\rho)$, we are automatically fixing a bundle $M\to M/G_0$, or
equivalently, a homotopy class in $[M/G_0, BG_0]$.

\begin{theorem} If the space $X$ is compact, then
\feq{union}{\mathrm{Bundle}(X,\mathrm{Aut}_G\(G_0\times \(F\otimes
V\)\))\approx \bigsqcup_{M\in \mathrm{Bundle}(X,G_0)}
\mathrm{Vect}_G(M,\rho).}
\end{theorem}

\proof By theorem 5, there is an inclusion \feq{inclusion}{
\bigcup_{M\in \mathrm{Bundle}(X,G_0)}
\mathrm{Vect}_G(M,\rho)\hookrightarrow\mathrm{Bundle}(X,\mathrm{Aut}_G\(G_0\times
\(F\otimes V\)\)).}

Now we will construct an inverse to the map (\ref{inclusion}), so the fact that the last union is disjoint will follow. Let
\ff{ \Psi_{\alpha\beta}:\;
\(U_\alpha\cap U_\beta\)\to \mathrm{Aut}_G\(G_0\times \(F\otimes V\)\)}
be the transition functions of a bundle $\xi\in\mathrm{Bundle}(X,\mathrm{Aut}_G\(G_0\times \(F\otimes V\)\))$.
By lemma 2, there is a continuous projection
of groups \linebreak $pr:\mathrm{Aut}_G\(G_0\times \(F\otimes V\)\)\to G_0$. So, by composition with
$pr$ we obtain a bundle with the discrete fiber $G_0$, and it is well known that $G_0$ acts
fiberwise and freely over the total space $M$ of this bundle and that $ M/G_0=X$.

 Also, we can assume that we have chosen an atlas such that there is a homeomorphism
 \ff{M\approx \underset{\alpha}{\bigcup}\(U_\alpha\times G_0\)\approx\underset{\alpha}{\bigcup}\(\bigsqcup_{[g]\in G_{0}}[g]U_{\alpha}\)}
 where the intersections are defined by the rule
 \ff{[1]U_\alpha\cap [g_{\alpha\beta}]U_\beta \approx U_\alpha\cap U_\beta}
 where $[g_{\alpha\beta}]=pr\circ \Psi_{\alpha\beta}$.

 On the other hand, we have
 \ff{\xi\approx \underset{\alpha}{\bigcup}\(U_\alpha\times \(G_0\times \(F\otimes V\)\)\)}
 where $U_\alpha\times \(G_0\times \(F\otimes V\)\)$ intersects $U_\beta\times \(G_0\times \(F\otimes V\)\)$
on the points $(x,g,f\otimes v)=(x,\Psi_{\alpha\beta}([g],f\otimes v))=(x,[g_{\alpha\beta}g],A_{\alpha\beta}[g](f\otimes v))$
where
\f{x\in U_\alpha\cap U_\beta}
and, once again,  we are using lemma 2 for the description of the operators $\Psi_{\alpha\beta}$.

Taking into account the homeomorphism
\ff{U_\alpha\times G_0\approx\bigsqcup_{[g]\in G_{0}}[g]U_{\alpha}}
we can rewrite
\ff{([g]x,f\otimes v)=([gg_{\alpha\beta}]x,A_{\alpha\beta}[g](f\otimes v))}.

 Therefore, the projection
 \ff{\(U_\alpha\times G_0\)\times \(F\otimes V\)\to U_\alpha\times G_0}
 extends to a well-defined and continuous projection
 \ff{\xi\to  M.}It is clear by the preceding formulas,
 that this projection will be $G$-equivariant, if $G$ acts canonically over
  the fibers and in by left translations on $G_0$ under the quotient map $G\to G/H=G_0$.
 So, we have $\xi\in \mathrm{Vect}_G(M,\rho)$.

 To end the proof, we make the remark that, by the theory of principal $G_0$-bundles, the construction
 of the space $M$ is up to equivariant homeomorphism. This means that the inverse to
 (\ref{inclusion}) is well defined. \qed


\begin{thebibliography}{20}
\bibitem{Mishchenko} Luke G., Mishchenko A. S., \textit{Vector Bundles And Their Applications.}
                                       Kluwer Academic Publishers Group (Netherlands), 1998.
                                       ISBN: 9780792351542
\bibitem{Conner} P. Conner, E. Floyd. \textit{Differentiable periodic maps.} Berlin, Springer-Verlag 1964.
\bibitem{Palais} Palais R.S. \textit{On the Existence of Slices for Actions of Non-Compact Lie Groups}
                            Ann.  Math., 2nd Ser., Vol. 73, No. 2. (1961), pp. 295-323.
\bibitem{Atiyah} Atiyah M.F., \textit{K-theory.} Benjamin, New York, (1967).
\bibitem{Serre}Serre J.P., \textit{Representations line\'{a}ires des groupes finis.} Hermann, Paris. 1967.
\bibitem{Levine} Levine M.,  Serp\'{e} C.,\textit{On a spectral sequence for equivariant K-theory}
                                                                        K-Theory (2008) 38 pp. 177–222

\end{thebibliography}
\end{document}